\documentclass[11pt]{article}

\usepackage{amsfonts}

\begin{document}

\baselineskip = 18 pt   

\parskip = \the\baselineskip


\parindent = 0 pt

\newtheorem{thm}{Theorem}[section]

\newtheorem{lem}[thm]{Lemma}
\newtheorem{definition}[thm]{Definition}

\newtheorem{pre-note}[thm]{Note}
\newenvironment{note}{\begin{pre-note}\rm}{\end{pre-note}}

\newtheorem{pre-proof}[thm]{Proof}
\newenvironment{proof}{\begin{pre-proof}\rm}{$\quad\bullet$\end{pre-proof}}

\newtheorem{pre-example}[thm]{Example}
\newenvironment{example}{\begin{pre-example}\rm}{\end{pre-example}}

\newtheorem{pre-numeq}[thm]{(}
\newenvironment{numeq}{\begin{pre-numeq}\rm\bf)$\quad\displaystyle}{$\end{pre-numeq}}


\newcommand{\RR}{\mathbb{R}}
\newcommand{\CC}{\mathbb{C}}
\newcommand{\HH}{\mathbb{H}}
\newcommand{\ZZ}{\mathbb{Z}}
\newcommand{\NN}{\mathbb{N}}

\newcommand{\Real}{\mathop{\rm Re}}
\newcommand{\Span}{\mathop{\rm Span}}
\newcommand{\inn}{\mathop{_{\rm \! inn}}}
\newcommand{\out}{\mathop{_{\rm \! out}}}
\newcommand{\naive}{\mathop{_{\rm \! naive}}}

\begin{flushright}
{\small
PUMA--563\\
$<$http://arXiv.org/abs/math/9911175$>$}
\end{flushright}

\begin{center}
{\Large SCATTERING FOR TIME SERIES WITH AN APPLICATION TO THE ZETA FUNCTION OF AN ALGEBRAIC CURVE}\\
\vskip 1 cm
{\Large Jean-Fran\c{c}ois Burnol}\\
\vskip 0.5 cm
November 1999
\end{center}

\vfill

I explain how the Lax-Phillips theory can be applied to a purely innovating time series and compute the corresponding scattering function. I then associate such a time series to an algebraic curve (of genus at least $1$) over a finite field and show that the Riemann Hypothesis (proven long ago) holds if and only if the scattering is causal (this causality is not independently established, though).
\vfill

\begin{flushleft}
{\small
Universit\'e de
Nice-Sophia-Antipolis\\
Laboratoire J.-A. Dieudonn\'e\\
Parc
Valrose\\
F-06108 Nice C\'edex 02\\
France\\
burnol@math.unice.fr\\}
\end{flushleft}

\clearpage


\section{Lax--Phillips scattering for innovating time series}

Let us first briefly review some concepts related to time series. See for example \cite{Grenander} for this as well as \cite{Dym-McKean} and \cite{Hoffman} for the underlying key results of classical harmonic analysis (Herglotz's theorem, Hardy spaces, Szeg\"o's theorem, inner and outer factors, the Beurling--Lax description of invariant subspaces, etc\dots) which we will use below.

Let $X = (X_n)_{n\in\ZZ}$ be a sequence of square-integrable random variables on some probability space, with zero mean, and such that the covariance ${\bf E}(\overline{X_n}X_m)$ depends only on the difference $n-m$ (and is denoted $\gamma_{n-m}$). We also assume that $\gamma_0 > 0$.

The $\gamma_j$'s satisfy the positivity condition of Herglotz and are thus the Fourier coefficients of a uniquely determined positive measure $\mu$ on the circle $S^1$: $\gamma_j = \int z^{-j} d\mu$ ($0 < \mu(S^1) < \infty$), the so-called ``spectral measure'' of the (weakly) stationary time series $X$.

Let $H(X)$ (the history of $X$) be the Hilbert space spanned by the $X_n$'s. There is a unique isometry $L^2(S^1,d\mu) \rightarrow H(X)$ which sends the function $z^n$ to the random variable $X_n$. Multiplication by $z$ then provides a (unique) isometry $U$ of $H(X)$ such that $U(X_n) = X_{n+1}$. We have a chain of increasing subspaces $H_n \subset H_{n+1}$ with $H_n$ spanned by the $X_m$'s, $m\leq n$. Let $H_{-\infty}$ be the intersection of all the $H_n$'s (the distant past of the process $X$). Then either $H_{-\infty} = H(X)$, or all inclusions $H_n \subset H_{n+1}$ are strict. If $H_{-\infty} = \{0\}$ we are necessarily in this second case and we will say that the time series $X$ is ``purely innovating'' (the traditional terminology is ``purely non-deterministic''). It then has no distant past, and no distant future either as $\alpha(z)\mapsto\overline{\alpha(z)}$ on $L^2(S^1,d\mu)$ gives an anti-unitary which exchanges $X_n$ with $X_{-n}$.

As an example of a stationary time series let us take a (variance $1$) white noise process $Y = (Y_n)_{n\in\ZZ}$ (by this it is just meant $\gamma_j(Y) = 0$ for $j \neq 0$, $\gamma_0 = 1$, or equivalently that the spectral measure is normalized Lebesgue measure), and also a sequence $(c_n)_{n\in\ZZ}$ in $l^2(\ZZ)$ and define $X_n$ as the (bi-sided) ``moving average'' $\sum_k c_k Y_{n-k}$. The spectral measure of $X$ is then $|\psi({1\over z})|^2 {d\theta\over2\pi}$ (where $\psi(z) = \sum_n c_n z^n$ in $L^2(S^1,{d\theta\over2\pi})$, $z=e^{i\theta}$). Any absolutely continuous measure can be written in this form, so any process with an absolutely continuous measure can be written as a (bi-sided) moving average. If the coefficients $c_n$ are such that $c_n = 0$ for $n >> 0$ or for $n << 0$, the moving average is said to be one-sided.

\begin{thm}[see {\cite[chap.2]{Grenander}}]
The following conditions on $X$ are equivalent:
\begin{enumerate}
        \item $X$ can be represented as a one-sided moving average (with respect to a white noise process).
        \item $X$ is purely innovating.
        \item The spectral measure of $X$ is absolutely continuous: $d\mu = f(\theta){d\theta\over2\pi}$ and furthermore $\int \log f(\theta)\ {d\theta\over2\pi} > -\infty$.
        \item The spectral measure of $X$ is absolutely continuous: $d\mu = f(\theta){d\theta\over2\pi}$ and furthermore there exists $\psi(z) = \sum_{n\geq0} c_n z^n$ in the Hardy Space $\HH^2$ such that $f(\theta)$ is almost everywhere equal to $|\psi(z)|^2$.
\end{enumerate}
\end{thm}

Let us suppose from now on that $X$ is such a purely innovating process. Let:
$${\cal D}_-^0 = \overline{\Span\{X_n, n< 0\}}$$
$${\cal D}_+^0 = \overline{\Span\{X_n, n\geq 0\}}$$
Then ${\cal D}_-^0$, ${\cal D}_+^0$, and the unitary shift operator $U$ satisfy the Lax-Phillips axioms for scattering (\cite{Lax-Phillips}):
$$j\leq 0 \Rightarrow U^j({\cal D}_-^0)\subset{\cal D}_-^0$$
$$\bigwedge U^j{\cal D}_-^0 = 0\quad\overline{\bigvee U^j{\cal D}_-^0} = H(X)$$
$$j\geq 0 \Rightarrow U^j({\cal D}_+^0)\subset{\cal D}_+^0$$
$$\bigwedge U^j{\cal D}_+^0 = 0\quad\overline{\bigvee U^j{\cal D}_+^0} = H(X)$$
except for the orthogonality axiom ${\cal D}_-^0 \perp {\cal D}_+^0$ which would be satisfied only for a white noise process. The orthogonality axiom guarantees the causality of the scattering matrix (more on this later). For this reason we will say that the choice $\{{\cal D}_-^0, {\cal D}_+^0\}$ defines the \emph{naive} scattering associated with $X$. 

\begin{definition}
The \emph{(dual) scattering} associated with the purely innovating process $X$ is given by the following choice of ``incoming'' and ``outgoing'' spaces:
$${\cal D}_- = \overline{\Span\{X_n, n\geq 0\}}^\perp = ({\cal D}_+^0)^\perp$$
$${\cal D}_+ = \overline{\Span\{X_n, n < 0\}}^\perp = ({\cal D}_-^0)^\perp$$
\end{definition}

We note that the orthogonality condition for $\{{\cal D}_-, {\cal D}_+\}$ boils down to $({\cal D}_-^0)^\perp \subset {\cal D}_+^0$ (rather than ${\cal D}_+^0 \subset ({\cal D}_-^0)^\perp $).

As the codimension of $H_n(X)$ in $H_{n+1}(X)$ is $1$, the Lax--Phillips theory teaches us that there exist (unique up to constants) isometric intertwiners 
$$\phi_- : H(X) \rightarrow L^2(S^1,{d\theta\over2\pi})$$
$$\phi_+ : H(X) \rightarrow L^2(S^1,{d\theta\over2\pi})$$
between the action of $U$ and multiplication by $z$ and such that 
$$\phi_-({\cal D}_-) = (\HH^2)^\perp$$
$$\phi_+({\cal D}_+) = \HH^2$$

The unitary $S = \phi_+ \cdot \phi_-^{-1}$ from $L^2(S^1,{d\theta\over2\pi})$ to itself is called the scattering matrix (in its spectral representation). As $S$ commutes with multiplication by $z$ it is multiplication with a measurable function of unit modulus $s(\theta)$, which we will call the scattering function.

\begin{thm}
Let $X$ be a purely innovating process, and $\psi(z)$ an element of the Hardy space $\HH^2$ such that the spectral measure of $X$ is normalized Lebesgue measure multiplied by $|\psi(z)|^2$. Let $\psi\out$ be the \emph{outer} part $\psi$ (see \cite{Dym-McKean}, \cite{Hoffman}). The dual scattering function associated with the process $X$ is given as 
$$s(\theta) = {\ \overline{\psi\out(e^{i\theta})}\ \over\ \psi\out(e^{i\theta})\ }$$
\end{thm}

\begin{proof}
We first determine the \emph{naive} scattering function associated with the pair $\{{\cal D}_-^0, {\cal D}_+^0\}$. Let $\phi_+^0$ be the corresponding intertwiner to the outgoing spectral representation. The map $L^2(d\mu) \rightarrow L^2({d\theta\over2\pi})$ given by $\alpha(z) \mapsto \psi(z)\alpha(z)$ is a unitary embedding (commuting with multiplication by $z$). It is in fact onto, as by a well-known result the zero set of $\psi(z)$ has Lebesgue measure $0$. It sends ${\cal D}_+^0$ to $\psi\HH^2$ which by a theorem of Beurling is also equal to $\psi\inn\HH^2$, with $\psi\inn$ the \emph{inner} part of $\psi$ (an inner function on the circle is a measurable function of modulus $1$ which is almost everywhere the non-tangential boundary value of an analytic function in the interior of the unit disc, itself bounded in modulus by $1$). Division by $\psi\inn(z)$ is a unitary, so the operator $\phi_+^0$ is given as $\alpha(z) \mapsto \psi\out(z)\alpha(z)$ from $L^2(d\mu)$ to $L^2({d\theta\over2\pi})$. In the same manner the ``incoming spectral representer'' $\phi_-^0$ is just multiplication by $\overline{\psi\out(z)}$ from $L^2(d\mu)$ to $L^2({d\theta\over2\pi})$. The naive scattering function is thus $\psi\out(e^{i\theta})\cdot\left(\overline{\psi\out(e^{i\theta})}\right)^{-1}$. And the looked-for $s(\theta)$ is its inverse. The conclusion follows
\end{proof}

\begin{definition}
The purely innovating process has \emph{causal scattering} if the associated dual scattering function is a causal function (equivalently if the naive scattering function is an inner function).
\end{definition}

\begin{note}
A causal function is an inner function with respect to the \emph{exterior} domain $|z| > 1$ (including $\infty$). Sometimes ``causal function'' refers to the values taken in that domain (which we still denote by $s(z)$), so we should perhaps say ``boundary value of a causal function''. This causality condition is equivalent to $s\naive(e^{i\theta})$ being an inner function (with respect to $|z| < 1$) as the following relations show:
\begin{eqnarray*}
|s(e^{i\theta})| &=& 1\\
s\naive(e^{i\theta})\cdot s(e^{i\theta}) &=& 1\\
(|z| < 1;\mbox{ causal case})\qquad\overline{s\naive(\overline{z})} &=& s({1\over z})\\
\end{eqnarray*}
\end{note}

\begin{note}
The outer part of $\psi$ can be expressed as (the boundary values of the exponential of) an integral (see \cite{Hoffman}) involving only the modulus of $\psi$, so that it is possible to express $s(\theta)$ directly in terms of the spectral measure $\mu$ (up to an arbitrary multiplicative constant of modulus $1$). We don't write up this formula as we will not need it.
\end{note}

\begin{note}
Perhaps our definition of causal scattering is a little too narrow and we should allow a pole at $z = \infty$ (equivalently at $z = 0$ for the naive scattering function). Replacing either the incoming or the outgoing subspace with a suitable shift would eliminate such a pole. The criterion of the next section would then also apply to genus $0$ curves.
\end{note}

\section{The congruence zeta--function as a time series}

Let $C$ be a smooth, geometrically irreducible, complete algebraic curve with field of constants the field with $q$ elements. Let $g$ be its genus, $h(C)$ its class number (the number of distinct divisor classes in each degree $d \in\ZZ$), and $Z(T)$ ($T = q^{-s}$) its zeta function (for all of this and more, see for example \cite{Moreno}). We will use ${\cal D}$ to denote an equivalence class of divisors. Two integers are associated with each such ${\cal D}$: its degree $d({\cal D})$ which belongs to $\ZZ$ and a dimension $l({\cal D})$ which belongs to $\NN$.

Most of what follows can be done also for $g = 0$ but behaves in the end slightly differently, so we will assume $g \geq 1$.

Let us now define a set of coefficients $e_m$ indexed by $m \in \ZZ$:
\begin{eqnarray*}
(|m|\geq g):\qquad e_m &=& -\, q^{-{|m|\over2}}\\
(|m|< g):\qquad e_m &=& {1\over h(C)}\left[\sum_{d({\cal D}) = m + g - 1} q^{l({\cal D}) - {m\over2}}\right]\ - \left[q^{m\over2} + q^{-{m\over2}}\right]\\
\end{eqnarray*}
and a time series $X = (X_n)_{n\in\ZZ}$ as the bi-sided moving average
$$X_n = \sum_m e_m Y_{n-m}$$
where $Y$ is a white noise process.

It is apparent that the Fourier series $\sum_m e_m z^m$ represents a rational function (with poles at ${1\over\sqrt{q}}$ and $\sqrt{q}$), so that its modulus has an integrable logarithm and $X$ is a \emph{purely innovating} process to which the theory described before applies.

\begin{thm}
The Riemann Hypothesis for the algebraic curve $C$ holds if and only if the associated process $X$ has causal scattering.
\end{thm}
\begin{proof}
First of all, let us recall the expression for the zeta function $Z(T)$ which exhibits it as a rational function (see \cite[chap. 3]{Moreno})
$$(q-1)Z(T) = \sum_{0\leq d({\cal D})\leq 2g-2} q^{l({\cal D})}T^{d({\cal D})} + h(C)\left[{q^g T^{2g-1}\over 1 -qT} - {1\over 1-T}\right]$$
and a straightforward calculation then shows the following identity in $L^2(S^1,{d\theta\over2\pi})$:
$$z^{-(g-1)}\,q^{g-1\over2}\,{q-1\over h(C)}\cdot Z({z\over\sqrt{q}}) = \sum_m e_m z^m$$
The spectral measure of $X$ is thus (note that $Z(\overline{T}) = \overline{Z(T)}$)
$$ d\mu = q^{g-1}\,\left({q-1\over h(C)}\right)^2\cdot |Z({z\over\sqrt{q}})|^2\ {d\theta\over2\pi}$$
It can be written  as $|\psi(z)|^2 {d\theta\over2\pi}$ for
$$\psi(z) = q^{g-1\over2}\,{q-1\over h(C)}\cdot {z - {1\over\sqrt{q}}\over1 - {z\over\sqrt{q}}}\cdot Z({z\over\sqrt{q}})$$
which belongs to $\HH^2$ as the only pole of $Z({z\over\sqrt{q}})$ in the open disc is at $z = {1\over\sqrt{q}}$. We now compute the naive scattering function as
$$s\naive(\theta) = {\ \psi\out(e^{i\theta})\ \over\ \overline{\psi\out(e^{i\theta})}\ }
= {\ \psi(e^{i\theta})\ \over\ \overline{\psi(e^{i\theta})}\ }
\cdot {\ \overline{\psi\inn(e^{i\theta})}\ \over\ \psi\inn(e^{i\theta})\ }
= {\psi(e^{i\theta})\over\psi(e^{-i\theta})}
\cdot \left(\psi\inn(e^{i\theta})\right)^{-2}$$
where $\overline{\psi(z)} = \psi(\overline{z})$ and $|\psi\inn(e^{i\theta})| = 1$ were used. The functional equation reads
$$Z({1\over qT}) = q^{1-g}\,T^{2-2g}\ Z(T)$$
which translates into
$$\psi({1\over z}) = \left({1-{z\over\sqrt{q}}\over z-{1\over\sqrt{q}}}\right)^2 z^{2-2g}\psi(z)$$
We obtain, for $|z| = 1$:
$$s\naive(z) = \left({z - {1\over\sqrt{q}}\over1 - {z\over\sqrt{q}}}\right)^2 z^{2g-2} \cdot \left(\psi\inn(z)\right)^{-2}$$
Let us assume that $s\naive(z)$ is in fact an inner function. The rational function on the right--hand side can be inner only if it has no poles in the open disc. But this means that $\psi\inn$ has no zeroes (as $0$ and $1\over\sqrt{q}$ are not acceptable candidates). This is just a way of phrasing the Riemann Hypothesis. And conversely under the Riemann Hypothesis the naive scattering function is
$$s\naive(z) = \left({z - {1\over\sqrt{q}}\over1 - {z\over\sqrt{q}}}\right)^2 z^{2g-2}$$
which is indeed an inner function
\end{proof}

Similar things can be done in the number field case too, but the Tate Gamma functions at infinite places are not of bounded characteristic in the half-plane $\Real(s) > {1\over2}$ and this is a source of additional hurdles.
In \cite{Burnol} a $p$--adic scattering problem was studied and the causality was established in that local setting.

\begin{flushleft}
{\small
Jean-Fran\c{c}ois Burnol\\
Universit\'e de Nice-Sophia-Antipolis\\
Laboratoire J.-A. Dieudonn\'e\\
Parc Valrose\\
F-06108 Nice C\'edex 02\\
France\\
burnol@math.unice.fr}\\
\end{flushleft}

\clearpage

\end{document}